\documentclass[12pt]{article}
\usepackage{amsmath,amssymb}
\def\la{\langle}
\def\ra{\rangle}
\makeatletter

\@addtoreset{equation}{section}
\makeatother
\begin{document}
\title{Ramanujan's $_1\psi_1$ summation theorem --- perspective, announcement of bilateral
$q$-Dixon--Anderson and $q$-Selberg integral extensions, and context\\
}
\author{
M{\small ASAHIKO} I{\small TO}%
\footnote{
School of Science and Technology for Future Life,
Tokyo Denki University,
Tokyo 120-8551, Japan 
\quad email: {\tt mito@cck.dendai.ac.jp}}
\ and P{\small ETER} J. F{\small ORRESTER}%
\footnote{
Department of Mathematics and Statistics, The University of Melbourne, Victoria 3010, Australia 
\quad email: {\tt p.forrester@ms.unimelb.edu.au}
}
}
\date{}
\maketitle
\begin{abstract}
\noindent 
The Ramanujan $_1 \psi_1$ summation theorem in studied from the perspective of $q$-Jackson integrals, $q$-difference equations and
connection formulas.  This is an approach which has previously been shown to  yield Bailey's very-well-poised $_6 \psi_6$ summation. Bilateral Jackson integral 
generalizations of the Dixon--Anderson
and Selberg integrals relating to the type $A$ root system are identified as natural candidates for multidimensional generalizations of the  Ramanujan $_1 \psi_1$ summation theorem.
New results of this type are announced, and furthermore they are put into context by reviewing from previous literature explicit product formulas for
Jackson integrals relating to other roots systems obtained from the same perspective.
\end{abstract}

{\scriptsize 2010 {\it Mathematical Subject Classification.} Primary 33D15, 33D67; Secondary 39A13}

{\scriptsize  {\it Keywords.} 
Dixon-Anderson integral; 
Selberg integral; 
Ramanujan's $_1\psi_1$ summation formula; 
Bailey's very-well-poised $_6\psi_6$ summation formula; 
}

\section{Introduction}
\label{S1}
%
%
Ramanujan's $_1\psi_1$ summation theorem is 
\begin{equation}
\label{eq:1psi1}
\sum_{\nu=-\infty}^\infty\frac{(a)_\nu}{(b)_\nu}\,x^\nu
=\frac{(ax)_\infty(q)_\infty(b/a)_\infty(q/ax)_\infty}
{(x)_\infty(b)_\infty(q/a)_\infty(b/ax)_\infty}.
\end{equation}
Here $(u)_\infty:=\prod^\infty_{l=0}(1-uq^l)$  
and with $0<q<1$,
$
(u)_\nu = (u)_\infty/(uq^\nu)_\infty, \nu \in \mathbb{Z}
$.
For generic parameters $a$, $b$, (\ref{eq:1psi1}) requires $|b/a|<|x|<1$ for absolute convergence.

As a part of a feature in the Notices of the AMS during December 2012, commemorating 125 years since Ramanujan's birth, S.O. Warnaar \cite{Wa12} recently reviewed aspects of (\ref{eq:1psi1}). One of these, important to us, is that the Ramanujan $_1\psi_1$ summation can be written as a Jackson integral ($q$-series), 
generalizing that corresponding to the $q$-binomial series.  Another, also important to us, is that (\ref{eq:1psi1}) permits multi-dimensional extensions. 
This in this article we aim to detail a perspective on (\ref{eq:1psi1}) in the setting of Jackson integrals, $q$-difference equations and connection formulas. 
We will then identify a natural setting, namely that of $q$-generalizations of
the Dixon--Anderson and Selberg integrals --- both multidimensional generalizations of the
Euler beta integral --- as candidates for permitting analysis from this perspective and as providing new multi-dimensional extensions of (\ref{eq:1psi1}).
We announce examples of such extensions, with the details being deferred to separate publications. Moreover, we relate these results to a broader
class of multidimensional Jackson integrals permitting analogous analysis, and relating to root systems.

One of the key features of our viewpoint --- that of establishing a $q$-difference equation for the left-hand side of (\ref{eq:1psi1}) ---
can be traced back to Mellin  in his study of the Gauss hypergeometric equation (see \cite[\S 1.5.3]{AK11}. However, it wasn't until the early 1990's in the work of Aomoto and Aomoto--Kato that a general theory relating solutions of the $q$-difference equation with special asymptotic behaviour, to a linear combination of a particular basis of solutions via a connection formula, was developed
\cite{Ao90}--\cite{AK11}.    
This has the consequence of providing an explanation for various special function formulae relating to $q$-hypergeometric functions such as (\ref{eq:1psi1}), and moreover of providing a method to explore the vast field of multi-dimensional generalizations. 

We begin in Section \ref{S2} by detailing this method as it applies to Ramanujan's $_1\psi_1$ summation theorem. We then demonstrate its generality by sketching its application to the derivation of  Bailey's very-well-poised $_6\psi_6$ summation. Moreover, it is also the case that various transformation formulae between $q$-hypergeometric series can be understood as connection formulae, and this is briefly mentioned. Section \ref{S3} 
introduces bilateral $q$-generalizations of the Dixon-Anderson and Selberg integrals relating to the type $A$ root system, and the application of our viewpoint to providing explicit product
evaluations is sketched. We put these new results in context by briefly reviewing explicit product  formulas for Jackson integrals relating
to other root systems, previously obtained by application of this method.

\section{Special function formulae for bilateral $q$-hypergeometric series}\label{S2}
\subsection{ Ramanujan's ${}_1 \psi_1$ summation theorem}\label{S2.1}
A textbook treatment of (\ref{eq:1psi1}) can be found in \cite[p.\,138]{GR04}. As noted in \cite{Wa12},
it is originally recorded (without proof) as Item 17 of Chapter 16 in the second of his three notebooks \cite[p.\,32]{Be91}, \cite{Ra57} although the notation is different. For us a notation which replaces the left-hand side in (\ref{eq:1psi1}) by a type of Jackson integral is most fundamental. 

For general $a \in \mathbb{C}$ the latter is defined by 
\begin{equation}
  \label{eq:1.2}
  \int_0^a f(z)d_qz=(1-q)\sum_{\nu=0}^\infty f(aq^\nu)aq^\nu
\end{equation}
Note that as $q \to 1^-$, 
$
\int_0^a f(z)d_qz \to \int_0^a f(z)dz
$. The Jackson integral (\ref{eq:1.2}) is formally extended to a bilateral series by defining
\begin{equation}
  \label{eq:1.3}
  \int_0^{a\infty} f(z)\frac{d_qz}{z}=(1-q)\sum_{\nu=-\infty}^\infty f(aq^\nu)
\end{equation}
when the right-hand side converges. In terms of the notation (\ref{eq:1.3}), it was noted by Askey \cite{As79} that with the replacement $(a,b,x)\mapsto (\xi q^\beta, \xi q,q^\alpha)$ on (\ref{eq:1psi1}), i.e.
\begin{equation}
  \label{eq:1.3a}
  I(\xi):=\int_0^{\xi\infty} z^\alpha\frac{(qz)_\infty}{(q^\beta z)_\infty}\frac{d_q z}{z}=C\ \frac{\xi^\alpha\theta(q^{\alpha+\beta}\xi)}{\theta(q^\beta\xi)}
\end{equation}
with $C:=(1-q)(q)_\infty(q^{1-\beta})_\infty/(q^\alpha)_\infty(q^{1-\alpha-\beta})_\infty$, and where $\theta (z) := (z)_\infty(q/z)_\infty(q)_\infty$.

The notation (\ref{eq:1.3a}) separates off the $\xi$-dependent factor. Its explicit form, as given on the right-hand side of (\ref{eq:1.3a}), is easy to verify. Thus taking into account the poles of $I(\xi)$, one sees that it can be expressed as 
$
\xi^\alpha f(\xi)/\theta (q^\beta\xi)
$
where $f(\xi)$ is some holomorphic function on $\mathbb{C}^*$. Since $I(\xi)$ is invariant under the shift $\xi \to q\xi$, the holomorphic function of $f(\xi)$ satisfies the $q$-difference equation 
$
f(q\xi)=-f(\xi)/q^{\alpha+\beta}\xi
$.
The solution of this equation is uniquely determined up to a constant (i.e.\,term independent of $\xi$) $C$ as 
$
f(\xi)=C\theta(q^{\alpha+\beta}\xi)
$, and the right-hand side of (\ref{eq:1.3a}) follows, but with $C$ still to be determined.

With the $\xi$-dependence of $I(\xi)$ known, to calculate $C$ we are free to choose a particular value of $\xi$. A convenient choice is $\xi=1$. Then the definition (\ref{eq:1.3}) shows that all the negative $\nu$ terms in the Jackson integral defining $I(\xi)$ vanish, and we have
\begin{equation}
  \label{eq:qB}
I(1) = \int_0^1z^\alpha\frac{(qz)_\infty}{(q^\beta z)_\infty}\frac{d_q z}{z}.  
\end{equation}
In fact (\ref{eq:qB}) is precisely the $q$-beta integral (see e.g. \cite[p.\,20]{GR04}), and the standard notation is $I(1) = B_q(\alpha,\beta)$. We are thus faced with the problem of evaluating the $q$-beta integral.

For this purpose we write the Jackson integral $I(\xi)$ as $I(\alpha;\xi)$ to emphasize the $\alpha$-dependence. Crucial to the evaluation of the $q$-beta integral is that $I(\alpha;\xi)$ satisfies the $\xi$-independent $q$-difference equation
\begin{equation}
  \label{eq:1.7}
  I(\alpha;\xi)=\frac{1-q^{\alpha+\beta}}{1-q^\alpha}I(\alpha+1;\xi).
\end{equation}
This $q$-difference equation can be derived from a viewpoint developed by Mellin in 1907 (see \cite[p.268]{AK11}).

Define the integrand in the definition (\ref{eq:1.3a}) by 
$\Phi(z):=z^\alpha(qz)_\infty/(q^\beta z)_\infty$. Then we see that
$$
b(z):=\frac{\Phi(qz)}{\Phi(z)}
= \frac{q^{\alpha}(1-q^\beta z)}{(1-qz)}
=\frac{b^+(z)}{b^-(qz)},
$$
where $b^+(z)=q^\alpha(1-q^\beta z)$ and $b^-(z)=1-z$. Next, for an arbitrary meromorphic function $\varphi(z)$ on $\mathbb{C}^*$, we define the symbol $\la\varphi(z)\ra$ by
$$
\la\varphi(z)\ra:
=\int_0^{\xi\infty}\varphi(z)\Phi(z)\frac{d_q z}{z}.
$$
One observes that  in this setting, and with the operator $\nabla$ defined by 
$\nabla\varphi(z):=\varphi(z)-b(z)\varphi(qz),$ the equation
\begin{equation}
\label{eq:nabla=0}
\la\nabla\varphi(z)\ra=0,
\end{equation}
holds true. This is because (\ref{eq:nabla=0}) is equivalent to the trivial equation 
$\la\varphi(qz)\ra=\la\varphi(z)\ra$, where it is assumed $\la\varphi(z)\ra$ converges. 
In particular, if we put $\varphi(z)=b^-(z)$, 
then $\nabla\varphi(z)=b^-(z)-b^+(z)
=(1-q^{\alpha+\beta})z-(1-q^\alpha)$ and thus (\ref{eq:nabla=0}) gives
$$
(1-q^{\alpha+\beta})\la z\ra-(1-q^\alpha)\la 1\ra=0.
$$
Noting from the definitions that $\la 1\ra=I(\alpha;\xi)$ and $\la z\ra=I(\alpha+1;\xi)$, the $q$-difference equation (\ref{eq:1.7}) follows.

The $q$-difference equation (\ref{eq:1.7}) is independent of $\xi$. The special feature of the case $\xi=1$, which follows from (\ref{eq:qB}), is that one has the explicit $\alpha \to \infty$ behaviour 
\begin{equation}
  \label{eq:1.9}
  I(\alpha;1)\sim (1-q)(q)_\infty/(q^\beta)_\infty
\end{equation}
(i.e.\,the $\nu = 0$ term in the definition (\ref{eq:1.2})). Equivalently with $\alpha$ fixed, $I(\alpha+N;1)$ has the asymptotic behaviour given by the right-hand side of (\ref{eq:1.9}). By iterating (\ref{eq:1.7}) it thus follows 
\begin{equation}
  \label{eq:1.9a}
I(\alpha;1)=\frac{(q^{\alpha+\beta})_N}{(q^\alpha)_N}I(\alpha+N;1)
=\frac{(q^{\alpha+\beta})_\infty}{(q^\alpha)_\infty}\lim_{N\to \infty}I(\alpha+N;1)
=(1-q)\frac{(q^{\alpha+\beta})_\infty(q)_\infty}{(q^\alpha)_\infty(q^\beta)_\infty}.
\end{equation}
Substituting (\ref{eq:1.9a}) in (\ref{eq:1.3a}) with $\xi=1$ gives the value of $C$ stated below the latter. 

We remark that (\ref{eq:1.3a}) can equivalently be written
$$
I(\xi)=\xi^\alpha\frac{\theta(q^{\alpha+\beta}\xi)\theta(q^\beta)}{\theta(q^{\alpha+\beta})\theta(q^\beta\xi)}I(1),
$$
This emphasizes the viewpoint of the above analysis as computing $I(\alpha;\xi)$ via a connection formula with the special solution $I(\alpha;1)$.
On a different front, we remark too that there has been recent application of (\ref{eq:1psi1}) in mathematical physics, in particular in the asymptotic
analysis of the fluctuation of the right most particle in the partially asymmetric exclusion process \cite{SS10}.

\subsection{Bailey's $_6\psi_6$ summation and some transformation formulae}
\label{S2.2}
In terms of notation $(a_1,\dots,a_r)_\nu:=(a_1)_\nu\cdots(a_r)_\nu$ the basic hypergeometric series $_r\psi_r$ is specified by 
\begin{equation}
  \label{eq:6.1}
  _r\psi_r
  \begin{bmatrix}
    \begin{matrix}
      a_1,\dots,a_r \\
      b_1,\dots,b_r
    \end{matrix} \,; q,x
  \end{bmatrix}
:= \sum_{\nu=-\infty}^\infty\frac{(a_1,\dots,a_r)_\nu}{(b_1,\dots,b_r)_\nu}x^\nu.
\end{equation}
With $0<q<1$ and generic parameters $\{a_i\},\{b_j\}$ this requires 
$|b_1\cdots b_r/a_1\cdots a_r|<|x|<1$ for absolute convergence. 
In terms of the notation (\ref{eq:6.1}) the left-hand side of (\ref{eq:1psi1}) is $_1\psi_1
\Big[\,
    \begin{matrix}
      a\\
      b
    \end{matrix}\, ; q,x
\Big]
$.

There is also a summation formula for the so-called very-well-poised case of $_6\psi_6$, evaluated at a particular value of $x$, due to Bailey (see \cite[p.\,140]{GR04}),
\begin{multline}
\label{eq:6.2}
_6\psi_6
\bigg[\,
    \begin{matrix}
      q\sqrt{a}, -q\sqrt{a}, b, c, d, e \\
      \sqrt{a}, -\sqrt{a}, aq/b, aq/c, aq/d, aq/e
    \end{matrix}\,; q,{\displaystyle\frac{a^2q}{bcde}}
\bigg]
\\ = \frac{(aq,aq/bc,aq/bd,aq/be,aq/cd,aq/ce,aq/de,q,q/a)_\infty}{(aq/b,aq/c,aq/d,aq/e,q/b,q/c,q/d,q/e,a^2q/bcde
)_\infty}.
\end{multline}
A derivation of (\ref{eq:6.2}) using $q$-difference equations, which parallels that for Ramanujan's $_1\psi_1$ summation presented above, can be given \cite{Ito06-2}.
The first point to note is that making the replacements \cite{IS08} 
\begin{equation}
  \label{eq:7.x}
  \sqrt{a} \mapsto \xi,\quad  (b,c,d,e)\mapsto (a_1\xi,a_2\xi,a_3\xi,a_4\xi)
\end{equation}
the left-hand side of (\ref{eq:6.2}) can be written as the Jackson integral 
\begin{equation}
  \label{eq:7.0}
 \frac{\xi^{\alpha_1+\alpha_2+\alpha_3+ \alpha_4 -1}}{(1-q)(1-\xi^2)} \prod_{i=1}^{4}\frac{(a_i\xi)_\infty}{(q\xi / a_i)_\infty}J(\xi), \quad
 J(\xi):= \int_0^{\xi\infty}\Phi(z)\Delta(z)\frac{d_qz}{z},
\end{equation}
where $a_i=q^{\alpha_i}$,  $\Delta(z):=z^{-1}-z$ and 
$$
\Phi(z):= \prod_{i=1}^4z^{\frac{1}{2}-\alpha_i}\frac{(qz/a_i)_\infty}{(za_i)_\infty}.
$$
Next, analogous to (\ref{eq:1.3a}), the dependence on $\xi$ can be determined to show \cite{Ito06-2}
\begin{equation}
  \label{eq:7.1}
  J(\xi)= \tilde{C}\frac{\xi\theta(\xi^2)}{\prod_{m=1}^4\xi^{\alpha_m}\theta(a_m\xi)}
\end{equation}
for some $\tilde{C}$ independent of $\xi$.

According to (\ref{eq:7.1}), to determine $\tilde{C}$ it suffices to evaluate $J(\xi)$ for a particular value of $\xi$. One sees that with $\xi=a_1$ the negative $\nu$ terms in the Jackson integral defining $J(\xi)$ vanishes and we have 
\begin{equation}
  \label{eq:7.2}
  J(a_1)=\int_0^{a_1}\Phi(z)\Delta(z)\frac{d_qz}{z}.
\end{equation}
In the simultaneous
limit 
\begin{equation}\label{8.0}
a_1\to a_1 q^{2N}, \quad a_2\to a_2 q^{-N},   \quad a_3 \to a_3 q^{-N},   \quad a_4 \to a_4 q^{-N} \: \: (N \to \infty)
\end{equation}
the asymptotic behaviour of $J(a_1)$ is simply given by the $\nu = 0$ term in (\ref{eq:7.2}). This is crucial since for general $\xi$ (see \cite[Corollary 6.2 with $n=1$]{Ito06-2} or \cite[Theorem 4.1 with $s=1$]{Ito09})
\begin{equation}
  \label{eq:8.1}
T_{a_i}J(\xi)=-
\frac{\prod_{k=1}^4(1-a_ia_k)}{a_i (1 - a_i^2) (1-a_1a_2a_3a_4)}J(\xi),
\end{equation}
where $T_{a_i}$ denotes the $q$-shift of the parameter $a_i$ such that $T_{a_i}:a_i\to qa_i$. Thus (\ref{eq:8.1}) can, in the case $\xi=a_1$, be solved for $J(a_1)$. 
By repeated use of (\ref{eq:8.1}) for the limit (\ref{8.0}), $J(a_1)$ is thus written as
$$
J(a_1) = (1 - q) {a_1^{1 - \alpha_1 -\alpha_2-\alpha_3- \alpha_4} (q)_\infty \prod_{2 \le i < j \le 4} (q a_i^{-1} a_j^{-1})_\infty \over
(q a_1^{-1} a_2^{-1} a_3^{-1} a_4^{-1})_\infty \prod_{k=2}^4 (a_1 a_k)_\infty},
$$
which is equivalent to Jackson's ${}_6 \phi_5$ summation formula \cite[p.44 (2.7.1)]{GR04}.
This allows $\tilde{C}$ in (\ref{eq:7.1}) to be determined 
as $\tilde{C}=(1-q)(q)_\infty\prod_{1 \le i < j \le 4} (q a_i^{-1} a_j^{-1})_\infty/(q a_1^{-1} a_2^{-1} a_3^{-1} a_4^{-1})_\infty$. 
Thus (\ref{eq:7.0}), which is the 
left-hand side of (\ref{eq:6.2}), has been evaluated in product form which upon recalling (\ref{eq:7.x}) is precisely the right-hand side of (\ref{eq:6.2}).

Not only can the $q$-difference equation method give a unified explanation of the basic hypergeometric summations (\ref{eq:1psi1}) and (\ref{eq:6.2}), we remark that various transformation formulae --- in particular those due to Sears \cite{Se51a, Se51b} and Slater \cite{Sl52,Sl66} --- can also be understood from the viewpoint \cite{IS08}. Briefly the transformation formulae are of the form of a basic bilateral series expanded as a linear combination of several specific bilateral series. Now the underlying difference equation is of rank higher than $1$. The transformation formulae are understood as connection formulae expressing a particular solution as a linear combination of solutions in a distinguished basis. 

\section{Higher dimensional generalizations}
\label{S3}
\subsection{Selberg and Dixon--Anderson integrals}
\label{S3.1}
A clue as to natural candidates for multi-dimensional generalizations comes from recalling that a special case of the $_1\psi_1$ summation is the $q$-beta integral (\ref{eq:qB}). In the limit $q\to 1$ this reduces to the Euler beta integral, while taking $q\to 1$ in the recurrence (\ref{eq:1.7}) implies the beta integral evaluation
\begin{equation}
  \label{eq:9.1}
  \int_0^1z^{\alpha-1}(1-z)^{\beta-1}dz = \frac{\Gamma(\alpha)\Gamma(\beta)}{\Gamma(\alpha+\beta)}.
\end{equation}
The significance of this for present purposes is that there are two multi-dimensional generalizations of (\ref{eq:9.1}): the Selberg  and Dixon--Anderson integrals. 

The former, discovered by Selberg in 1941 (see \cite{FW08} for historical aspects) reads
\begin{multline}
  \label{eq:2.1x}
\int_0^1dz_1 \ldots \int_0^1dz_n 
\prod_{i=1}^n z_i^{\alpha-1}(1-z_i)^{\beta - 1} 
\prod_{1\leq j<k \leq n}|z_k-z_j|^{2\tau}
\\
=\prod_{j=1}^n
\frac{\Gamma(\tau j +1)\Gamma(\alpha + (n-j)\tau)\Gamma(\beta+(n-j)\tau)}{\Gamma(\tau+1)\Gamma(\alpha +\beta + (n+j-2)\tau)},
\end{multline}
while the latter, due independently to Dixon \cite{Di05} and Anderson \cite{An91}  (the paper \cite{Di05}, written in 1905, only became widely known in recent times; see \cite{FW08}) reads
\begin{multline}
  \label{eq:2.2y}
  \int_{x_{n-1}}^{x_n}dz_n \ldots \int_{x_1}^{x_2}dz_2 \int_{x_0}^{x_1}dz_1\prod_{i=1}^n\prod_{j=0}^n|z_i-x_j|^{s_j-1}\prod_{1\leq k<l \leq n}(z_l-z_k)
\\
=\frac{\Gamma(s_0)\Gamma(s_1)\ldots\Gamma(s_n)}{\Gamma(s_0+s_1+\ldots+s_n)}\prod_{0\leq i<j\leq n}(x_j-x_i)^{s_i+s_j-1}.
\end{multline}
In the case $n=1$ both (\ref{eq:2.1x}) and (\ref{eq:2.2y}) reduce to the Euler beta integral (\ref{eq:9.1}).

Building on a number of earlier works \cite{As80,Ev92,Ev94,Ha88,Kad88}, the natural multi-dimensional Jackson integral generalizations of (\ref{eq:2.1x}) and (\ref{eq:2.2y}) which contains the Ramanujan $_1\psi_1$ summation when $n=1$, has recently been identified to be given by \cite{IF13}
\begin{equation}
  \label{eq:2.4z}
\tilde{I}(\xi):= \int_0^{\xi\infty}(z_1 \cdots z_n)^\alpha \prod_{i=1}^n \prod_{j=1}^m \frac{(qa_j^{-1}z_i)_\infty}{(b_jz_i)_\infty}\prod_{1 \leq k < l \leq n}z_k^{2\tau-1}\frac{(q^{1-\tau}z_l/z_k)_\infty}{(q^\tau z_l/z_k)_\infty}(z_k-z_l)\frac{d_qz_1}{z_1}\land \cdots \land \frac{d_qz_n}{z_n}.
\end{equation}
Here $\alpha, \tau \in \mathbb{C}$, $a_1, \ldots, a_m, b_1, \ldots, b_m \in \mathbb{C}^*$, $\xi = (\xi_1, \ldots, \xi_n) \in (\mathbb{C}^*)^n$ and 
$$\int_0^{\xi\infty}f(z)\frac{d_qz_1}{z_1}\land \cdots \land \frac{d_qz_n}{z_n} 
:= (1-q)^n \sum_{(\nu_1,\ldots,\nu_n)\in \mathbb{Z}^n}f(\xi_1q^{\nu_1},\ldots,\xi_nq^{\nu_n}).$$
\par
With $\tau$ a positive integer, $m=1$ and $\xi= (1, q^\tau, \ldots, q^{(n-1)\tau})$, taking the limit $q \to 1$ reclaims the Selberg integral (\ref{eq:2.1x}). This same setting, but with $m=2$  and $\alpha = 1$ also reclaims (\ref{eq:2.1x}) in the limit $q \to 1$, although with a change of variables such that the terminals and integrand is translated from $[0,1]^n$ to $[a_1,a_2]^n$. 

With $\tau=1/2, m=n$, $\xi=(a_1, \ldots, a_n)$, $(\alpha, a_j, b_j)\mapsto (s_0,x_{j-1}, q^{j-1}/x_{j-1})$ taking the limit $q \to 1$ reclaims the Dixon--Anderson integral (\ref{eq:2.2y}) with $x_0=0$.
The case $\tau=1/2, m=n+1$, $\alpha = 1$ and other parameters as stated reduces to (\ref{eq:2.2y}) without the restriction $x_0 = 0$.

\subsection{Bilateral $q$-generalizations of the Selberg and Dixon--Anderson integrals}
\label{S3.2}
The strategy given in Section \ref{S2.1} for the derivation of (\ref{eq:1psi1}), involving first determining the $\xi$-dependence in (\ref{eq:1.3a}), deriving a $q$-difference equation in $\alpha$ for $I(\xi)$, and solving the difference equation for a particular $\xi$ by making use of the knowledge of the $\alpha \to \infty$ asymptotic behaviour of $I(\xi)$ for this value of $\xi$ again hold us in good stead for the evaluation of $\tilde{I}(\xi)$ in the cases it admits a limit to the Selberg or Dixon--Anderson integrals.

The $m=1$ case of (\ref{eq:2.4z}) is due to Aomoto \cite{Ao98}, who obtained the evaluation
\begin{equation}
  \label{eq:2.5}
  \tilde{I}(\xi) = c_0 \prod_{i=1}^n \xi_i^{\alpha + 2(n-i)\tau}\frac{\theta(q^{\alpha +(n-1)\tau}b_1\xi_i)}{\theta(b_1\xi_i)}\prod_{1\leq j<k \leq n }\frac{\theta(\xi_k/\xi_j)}{\theta(q^\tau \xi_k/\xi_j)},
\end{equation}
where $c_0$ is independent of $\xi$ and is explicitly given by 
$$
c_0 = \prod_{j=1}^n \frac{(1-q)(q)_\infty(q^{1-j\tau})_\infty(q^{1-(j-1)\tau}a_1^{-1}b_1^{-1})_\infty}{(q^{1-\tau})_\infty(q^{\alpha+(j-1)\tau})_\infty(q^{1-\alpha-(n+j-2)\tau}a_1^{-1}b_1^{-1})_\infty}.
$$
The case $m=2$ and $\alpha=1$, which like the $m=1$ case also admits a limit to the Selberg integral (\ref{eq:2.1x}) as remarked above, is for general $\tau$ a great deal more complicated. A product formula analogous to the right-hand side of (\ref{eq:2.5}) only results by considering a sum of $n+1$ Jackson integrals (\ref{eq:2.4z}), each defined by a different $\xi$,
$$
\xi= (\underbrace{x_1,x_1q^\tau,\ldots,x_1q^{(i-1)\tau}}_i,\underbrace{x_2,x_2 q^\tau,\ldots,x_2 q^{(n-i-1)\tau}}_{n-i})
\in (\mathbb{C}^*)^n
\qquad(0\le i\le n)
$$
with $x_1, x_2 \in \mathbb{C}^*$, and furthermore weighted by certain functions of $x_1,x_2$ \cite{IF13}. In the special case $(x_1,x_2)=(a_1,a_2)$ and $\tau$ a positive integer this is shown in \cite{IF13} to be equivalent to a $q$-generalization of the Selberg integral conjectured by Askey \cite{As80} and subsequently proved by Evans \cite{Ev92}.

In relation to bilateral $q$-generalizations of the Dixon--Anderson integral, for $m=n, \tau=1/2$, $I(\xi)$ can be expressed as the ratio of theta functions
\begin{equation}
  \label{eq:2.6}
  \tilde{I}(\xi) = c_1 (\xi_1\cdots \xi_n)^\alpha \frac{\theta(q^\alpha\xi_1\cdots \xi_n b_1\cdots b_n)}{\prod_{i=1}^n \prod_{j=1}^n \theta(\xi_i b_j)}\prod_{1 \leq i < j \leq n}\xi_j \theta(\xi_i/\xi_j),
\end{equation}
where  $c_1$ is independent of $\xi$ and has the explicit evaluation
$$
c_1 = \frac{(1-q)^n(q)_\infty^n \prod_{i=1}^n \prod_{j=1}^n(qa_i^{-1}b_j^{-1})_\infty}{(q^\alpha)_\infty(q^{1-\alpha}a_1^{-1} \cdots a_n^{-1}b_1^{-1}\cdots b_n^{-1})_\infty}.
$$
We remark that (\ref{eq:2.6}) is equivalent to the Milne--Gustafson summation formula \cite{Gu87,Mi86}. The second case of (\ref{eq:2.4z}) permitting a limit to the Dixon--Anderson integral, namely $m=n+1,\alpha=1$, is like the second of (\ref{eq:2.4z}) permitting a limit to the Selberg integral reviewed above, more complicated in its structure with a product formula resulting upon summing $n+1$ Jackson integrals (\ref{eq:2.4z}), weighted by $(-1)^{i-1}$ and with $\xi$ given by 
$$
\xi=(x_1, \ldots, x_{i-1}, x_{i+1}, \ldots, x_{n+1})\in (\mathbb{C}^*)^n
\qquad(1\le i\le n+1).
$$
In the special case $a_i = x_i$ $(i=1, \ldots, n+1)$ this is shown in \cite{IF12} to be equivalent to a $q$-generalizations of the Dixon--Anderson integral due to Evans \cite{Ev94}.

\subsection{Multi-dimensional Jackson integrals associated with root systems}
\label{S3.3}

We have seen that the $_6\psi_6$ summation in Bailey's formula (\ref{eq:6.2}) can be written in terms of the Jackson integral in (\ref{eq:7.0}). As with the Jackson integral in (\ref{eq:1.3a}),  it is possible to generalize (\ref{eq:7.0}) to a bilateral multi-dimensional Jackson integral which permits analysis along the lines of the strategy given in Section \ref{S2.1}. These are all special cases of the family of Jackson integrals  
\begin{equation}
  \label{eq:13.1}
   \tilde{J}(\xi)= \int_0^{\xi\infty}\tilde{\Phi}(z)\tilde{\Delta}(z)\frac{d_qz_1}{z_1}\land \cdots \land \frac{d_qz_n}{z_n},
\end{equation}
where 
$$
\tilde{\Phi}(z) = \prod_{i=1}^n \prod_{m=1}^{2s+2}z_i^{1/2 - \alpha_m} \frac{(qa_m^{-1}z_i)_\infty}{(a_mz_i)_\infty} 
\prod_{1 \leq j < k \leq n}z_j^{1-2\tau} \frac{(qt^{-1}z_j/z_k)_\infty}{(tz_j/z_k)_\infty} \frac{(qt^{-1}z_jz_k)_\infty}{(tz_jz_k)_\infty}, 
$$
$$
\tilde{\Delta}(z) = \prod_{i=1}^n \frac{1-z_i^2}{z_i} \prod_{1 \leq j < k \leq n} \frac{(1-z_j/z_k)(1-z_jz_k)}{z_j}
$$
with $a_m=q^{\alpha_m}$ and $t=q^\tau$. 
In particular the case $s=1$ (general $n \in \mathbb{Z}^+)$ permits a product function evaluation \cite{Ito06-2} 
equivalent to earlier results due to Gustafson \cite{Gu90} and van Diejen \cite{vD97}. 
The case $n=1$ (general $s \in \mathbb{Z}^+$) permits a connection formula \cite{Ito08} 
generalizing the Sears--Slater transformation formula  
as mentioned at the end of Section \ref{S2.2}.

The symmetries of the summand in (\ref{eq:13.1}) contrast with those of (\ref{eq:2.4z}). The latter relate to the $A$-type root system, while the former relate to the $BC$-type root system. 
In fact it is similarly true that bilateral Jackson integrals with summands corresponding to 
any irreducible reduced root systems, and which permit product formula evaluations, can also be formulated \cite{Ito01,Ito03,Ito06-1,IT08,IT10}.

\subsection*{Acknowledgements}
We thank Ron Tidhar for his help with the preparation of the manuscript. This work was supported by the Australian Research Council and JSPS KAKENHI Grant Number 25400118.

\end{document}